\documentclass[12pt,leqno]{article}
\usepackage{amsfonts}
\usepackage{amsmath, amsthm, amsfonts, amssymb, color}
\usepackage{mathrsfs}
\usepackage{color}
\theoremstyle{definition}

\newcommand{\scr}[1]{\mathscr #1}
\definecolor{wco}{rgb}{0.5,0.2,0.3}

\numberwithin{equation}{section} \theoremstyle{remark}

\newcommand{\ua}{\uparrow}

\title{{\bf Bi-Coupling Method and Applications }\footnote{Supported in
 part by  the National Key R\&D Program of China (No. 2022YFA1006000, 2020YFA0712900) and NNSFC (11921001).} }
\author{
{\bf Panpan Ren$^{b)}$,    Feng-Yu Wang$^{a)}$  }\\
\footnotesize{$^{a)}$ Center for Applied Mathematics, Tianjin University, Tianjin 300072, China}\\
 \footnotesize{$^{b)}$ Department of Mathematics, City University of  Hong Kong, Tat Chee Avenue, Hong Kong,  China }\\
\footnotesize{    wangfy@tju.edu.cn}}

\begin{document}
\allowdisplaybreaks
\def\R{\mathbb R}  \def\ff{\frac} \def\ss{\sqrt} \def\B{\mathbf
B}
\def\N{\mathbb N} \def\kk{\kappa} \def\m{{\bf m}}
\def\ee{\varepsilon}\def\ddd{D^*}
\def\dd{\delta} \def\DD{\Delta} \def\vv{\varepsilon} \def\rr{\rho}
\def\<{\langle} \def\>{\rangle}
  \def\nn{\nabla} \def\pp{\partial} \def\E{\mathbb E}
\def\d{\text{\rm{d}}} \def\bb{\beta} \def\aa{\alpha} \def\D{\scr D}
  \def\si{\sigma} \def\ess{\text{\rm{ess}}}\def\s{{\bf s}}
\def\beg{\begin} \def\beq{\begin{equation}}  \def\F{\scr F}
\def\Ric{\mathcal Ric} \def\Hess{\text{\rm{Hess}}}
\def\e{\text{\rm{e}}} \def\ua{\underline a} \def\OO{\Omega}  \def\oo{\omega}
 \def\tt{\tilde}\def\[{\lfloor} \def\]{\rfloor}
\def\cut{\text{\rm{cut}}} \def\P{\mathbb P} \def\ifn{I_n(f^{\bigotimes n})}
\def\C{\scr C}      \def\aaa{\mathbf{r}}     \def\r{r}
\def\gap{\text{\rm{gap}}} \def\prr{\pi_{{\bf m},\varrho}}  \def\r{\mathbf r}
\def\Z{\mathbb Z} \def\vrr{\varrho} \def\ll{\lambda}
\def\L{\scr L}\def\Tt{\tt} \def\TT{\tt}\def\II{\mathbb I}
\def\i{{\rm in}}\def\Sect{{\rm Sect}}  \def\H{\mathbb H}
\def\M{\mathbb M}\def\Q{\mathbb Q} \def\texto{\text{o}} \def\LL{\Lambda}
\def\Rank{{\rm Rank}} \def\B{\scr B} \def\i{{\rm i}} \def\HR{\hat{\R}^d}
\def\to{\rightarrow}\def\l{\ell}\def\iint{\int}\def\gg{\gamma}
\def\EE{\scr E} \def\W{\mathbb W}
\def\A{\scr A} \def\Lip{{\rm Lip}}\def\S{\mathbb S}
\def\BB{\scr B}\def\Ent{{\rm Ent}} \def\i{{\rm i}}\def\itparallel{{\it\parallel}}
\def\g{{\mathbf g}}\def\Sect{{\mathcal Sec}}\def\T{\mathbb T}\def\BB{{\bf B}}
\def\f{\mathbf f} \def\g{\mathbf g}\def\BL{{\bf L}}  \def\BG{{\mathbb G}}
\def\Bd{{D^E}} \def\BdP{D^E_\phi} \def\Bdd{{\bf \dd}} \def\Bs{{\bf s}} \def\GA{\scr A}
\def\Bg{{\bf g}}  \def\Bdd{\psi_B} \def\supp{{\rm supp}}\def\div{{\rm div}}
\def\ddiv{{\rm div}}\def\osc{{\bf osc}}\def\1{{\bf 1}}\def\BD{\mathbb D}\def\GG{\Gamma}
\def\H{{\bf H}}
\maketitle

\begin{abstract} By developing a new technique called  the bi-coupling argument, we estimate the relative entropy between different diffusion processes in terms of  the distances of initial  distributions and  drift-diffusion coefficients. As an application, the entropy-cost inequality is established  for
McKean-Vlasov SDEs with spatial-distribution dependent noise, which is open for a long time and has  potential applications in optimal transport, information theory and mean field particle systems.
  \end{abstract} \noindent

 \vskip 2cm

 \section{Introduction}

The main purpose of this paper is to establish the entropy-cost inequality for McKean-Vlasov SDEs with spatial-distribution dependent noise, which has been open for a long time due to the essential difficulty caused by
the distribution dependence of noise.
 To overcome this difficulty, we develop a new coupling argument, called bi-coupling, to cancel the short time singularity in the entropy upper bound for two diffusions presented in \cite{BRS}.

  In this part, we first introduce the background of the study  from applied areas including the information theory, optimal transport and mean field particle systems, then explain the main difficulty of the study, and finally figure out the main idea of the present study and the structure of the paper.

\subsection{Background of the study}

Let $\scr P_2$ be the space of probability measures on $\R^d$ having finite second moments,  which is a Polish space under the quadratic Wasserstein distance
\beq\label{AS0}\W_2(\mu,\nu):=\inf_{\pi\in\C(\mu,\nu)}\bigg(\int_{\R^d\times\R^d}|x-y|^2\pi(\d x,\d y)\bigg)^{\ff 1 2},
\ \ \mu,\nu\in\scr P_2,\end{equation}
where $\C(\mu,\nu)$ is the set of all couplings of $\mu$ and $\nu$. In the theory of optimal transport, $\W_2$ refers to the  optimal transportation cost induced by the quadratic cost function, and $(\scr P_2,\W_2)$ is called the Wasserstein space where a nice  analysis and geometry structure has been developed, see for instance Otto's celebrated paper \cite{OTTO} and Villani's  monograph \cite{V}.

In the information theory, the relative entropy functional describes the chaos of a distribution with respect to a reference measure, which refers to  the difference of  Shannon entropies for two distributions,
and  is known as  the Kullback-Leibler divergence or the information divergence \cite{KL}.
For two probability measures $\mu$ and $\nu$,  the relative entropy  of $\nu$ with respect $\mu$ is defined as
 $$\Ent(\nu|\mu):=\beg{cases} \int_{\R^d} \{\log\ff{\d\nu}{\d\mu}\}\d\nu,\ &\text{if} \ \ff{\d\nu}{\d\mu}\ \text{exists},  \\
 \infty,\ &\text{otherwise.} \end{cases} $$

Both Wasserstein distance and relative entropy have wide applications in applied areas including deep learning and Bayesian statistics \cite{GBC}.  When $\mu$ is the standard Gaussian measure on $\R^d$,   Talagrand   \cite{T} found the beautiful inequality
 $$\W_2(\nu,\mu)^2\le 2 \Ent(\nu|\mu),\ \ \nu\in \scr P_2,$$
 where the constant $2$ is sharp. This inequality was then extended in \cite{OV, BGL} as
 \beq\label{T} \W_2(\nu,\mu)^2\le C\, \Ent(\nu|\mu),\ \ \nu\in \scr P_2\end{equation}
  for a constant $C>0$ and
  a probability measure $\mu$ satisfying the log-Sobolev inequality
 $$\mu(f^2\log f^2):=\int_{\R^d}f^2\log f^2\d\mu\le C \mu(|\nn f|^2),\ \ f\in C_b^1(\R^d),\mu(f^2)=1.$$
 The inequality \eqref{T} enables one to estimate the cost from above by using the entropy.

 However, comparing with the Wasserstein distance $\W_2$,  the entropy is usually harder to estimate from above. For instance, $\W_2$ between the distributions of two SDEs can be bounded by the expectation of the distance square of the solutions, which is easily derived using It\^o's formula. But the entropy between solutions of SDEs is harder to estimate from above, since the heat kernels (distribution densities) are unknown. So, it is crucial to establish the inverse Talagrand inequality by bounding the entropy using  $\W_2$.

 In general,  $\W_2(\mu,\nu)^2$ can not dominate   $\Ent(\nu|\mu)$ since the former is finite for any $\mu,\nu\in \scr P_2$ but the latter becomes infinite  when $\nu$ is not absolutely continuous with respect to $\mu$.
So, to derive an inverse Talagrand inequality, we consider the entropy between two stochastic systems for which the entropy decays in time according to the H-theorem in information theory. In this spirit, a sharp   entropy-cost inequality  was found by the second named author \cite{W10} for diffusion processes on a manifold $M$. According to \cite[Theorem 1.1]{W10},
 for any constant $K\in\R$,  the  Bakry-Emery curvature of the diffusion process is bounded below by  $K\in \R$
 if and only if  the following entropy-cost inequality holds:
\beq\label{AS}\Ent(P_t^*\mu|P_t^*\nu)\le \ff K{2(\e^{2Kt}-1)}  \W_2^\rr(\mu,\nu)^2,\ \ \mu,\nu\in \scr P_2(M),\ \ t>0,\end{equation}
   where $P_t^*\mu$ is the distribution of the diffusion process at time $t$ with initial distribution $\mu$,
     $\W_2^\rr$ is the quadratic Wasserstein distance induced by the Riemannian distance $\rr$ on $M$
   (when $M=\R^d$ it reduces to $\W_2$ defined in \eqref{AS0}), and $\scr P_2(M)$ is the set of probability measures on $M$ having finite second moment.
 This inequality has applications for both short and long times:
 \beg{enumerate} \item[$\bullet$] For small time, \eqref{AS} describes an instant finite property of the entropy, i.e. even  though the initial entropy is infinite, the entropy at any time $t>0$ becomes finite, and the short time behavior of the entropy behaves like $t^{-1}$;
 \item[$\bullet$] For long time, \eqref{AS}   provides exponential decay of $P_t^*$  in entropy by using that in $\W_2$  which is easier to verify in applications. \end{enumerate}
The inequality \eqref{AS} is  equivalent to the  log-Harnack inequality   (see \cite{Wbook})
 $$P_t \log f(x)\le \log P_tf(y)+ \ff {K\rr(x,y)^2}{2(\e^{2Kt}-1)},\ \ x,y\in\R^d, t>0, f\in \B^+(M),$$
 where  $\B^+(M)$ is the space of all uniformly positive measurable functions  on $M$, and $P_tf(x):=\int_{\R^d} f(y) \d (P_t^*\dd_x)$ is the associated diffusion semigroup.
 As a member in the family of dimension-free Hananck inequalities (see \cite{W97,W10,W11}), the log-Harnack inequality has crucial applications in optimal transport, curvature  on Riemennian manifolds or  metric measure spaces,  see for instance \cite{AGS,RW, W10,Wbook}.

 In this paper, we aim to establish the entropy-cost inequality of type \eqref{AS} for the nonlinear Fokker-Planck equation on  $\scr P_2$:
\beq\label{FP} \pp_t \mu_t= (L_{t,\mu_t}^{a,b})^*\mu_t,\ \ t\in [0,T],\end{equation}
where $T>0$ is a fixed time, and for any $(t,\mu)\in [0,T]\times \scr P_2,$ $(L_{t,\mu}^{a,b})^*$ is the $L^2(\R^d)$-adjoint operator of
\beq\label{L0}L_{t,\mu}^{a,b} := \sum_{i,j=1}^d a^{ij} (t,\cdot,\mu)\pp_i\pp_j +\sum_{i=1}^d b^i(t,\cdot,\mu)\pp_i.\end{equation}   Recall that a  continuous map $\mu_\cdot: [0,T]\to\scr P_2$ is called a solution to \eqref{FP}, if
for any $f\in C_0^\infty(\R^d)$ we have  $\int_0^t |\mu_s(L_{s,\mu_s}f)|\d s<\infty$ and
$$\mu_t(f)=\mu_0(f) +\int_0^t \mu_s(L_{s,\mu_s}f)\d s,\ \ t\in [0,T].$$

By the propagation of chaos, see \cite{SN}, under   reasonable conditions we have
$$\mu_t=\lim_{N\to\infty} \ff 1 N \sum_{i=1}^N \dd_{X_t^{i,N}}\ \text{in}\ L^2(\OO\to \scr P_2;\P),$$
  where for every $N\in\mathbb N$, $(X_t^{i,N})_{1\le i\le N}$ is  the associated mean field particle system with $N$ many particles, and
$\mu_t$ is the distribution of the solution $X_t$ to the following McKean-Vlasov SDE:
\beq\label{E1} \d X_t= b(t,X_t, \L_{X_t})\d t+  \si(t,X_t, \L_{X_t})\d W_t,\ \ t\in [0,T], \end{equation}
where   $\L_{X_t}$ is the distribution of $X_t,$  $W_t$ is the $d$-dimensional Brownian motion under a standard probability base $(\OO,\F,\{\F_t\}_{t\in [0,T]},\P)$,  $\si:=\ss{2a}$, and $(a,b)$  comes  from  $L_{t,\mu}^{a,b}$ in \eqref{L0}.
According to \cite{BR}, under a mild integrability condition, \eqref{FP} is well-posed in $\scr P_2$ if and only if \eqref{E1} has weak well-posedness for distributions in $\scr P_2,$ and
in this case
$\mu_t=P_t^*\mu:=\L_{X_t^\mu} $  is the unique solution to the nonlinear Fokker-Planck equation \eqref{FP} with $\mu_0=\mu$, where $X_t^\mu$ solves \eqref{E1}   with $\L_{X_0}=\mu$.

We intend to find a constant $c>0$ such that
\beq\label{EC} \Ent(P_t^*\mu|P_t^*\nu)\le \ff c t \W_2(\mu,\nu)^2,\ \ t\in (0,T], \mu,\nu\in \scr P_2.\end{equation}
When the noise is distribution-free, i.e. $\si(t,x, \mu)=\si(t,x)$ does not depend on $\mu$,
\eqref{EC}  has already been derived and applied in the literature,  \eqref{EC} has been established in  \cite{HW18, HW22a, RW, FYW1, FYW3}  under different conditions,   see also \cite{HRW19,HS,WZK} for extensions to the infinite-dimensional and    reflecting models.
When the noise coefficient is also distribution dependent, the coupling by change of measures applied in the above references does not apply. Recently, for $\si(t,x,\mu)=\si(t,\mu)$ independent of the spatial variable $x$,  \eqref{EC} has been established in \cite{HW22} by using   a noise decomposition argument,  see also \cite{BH} for the study on a special model.

However, when the noise is spatial-distribution dependent,
this type inequality has been open  for a long time until the new coupling technique  (bi-coupling)  is developed in the present paper, for which we construct a new diffusion process which is coupled with the other two processes respectively, see Section 2 below for details.

 We would like to indicated that after an earlier version of this paper is available online (arXiv:2302.13500),  the bi-coupling method has been applied in \cite{QRW, HRW} to different models to derive new  estimates on entropy and probability distances,  so that the efficient and originality  of this new method has been illustrated.

\subsection{Existing entropy inequality and difficulty of the present study}
Noting that $P_t^*\mu$ is the distribution of the diffusion process $X_t^\mu$ generated by
$(L_{t, P_t^*\mu}^{a,b})_{t\in [0,T]}$, the left hand side in \eqref{EC}  is the entropy between the distributions of two diffusion processes  generated by
$L_{t,P_t^*\mu}^{a,b}$ and $L_{t,P_t^*\nu}^{a,b}$ respectively. So, the study reduces to estimate the entropy  between two different diffusion processes.

 In general,  let $\GG$ be the space of $(a,b)$, where
 $$b: [0,T]\times\R^d\to\R^d,\ \ \ a: [0,T]\times \R^d\to \R^d\otimes\R^d$$
 are measurable, and for any $(t,x)\in [0,T]\times \R^d,$ $a(t,x)$ is positive definite. For any $(a,b)\in \GG,$ consider the  time dependent second order differential operators on $\R^d$:
 $$L_{t}^{a,b}:=  {\rm tr}\{a(t,\cdot) \nn^2\}+ b(t,\cdot)\cdot\nn,\ \ t\in [0,T].$$
Let $(a_i,b_i)\in \GG, i=1,2,$ such that for any $s\in [0,T)$, each $(L_t^{a_i,b_i})_{t\in [s,T]}$ generates a unique diffusion process $(X_{s,t}^{i,x})_{(t,x)\in [s,T]\times\R^d}$ with $X_{s,s}^{i,x}=x$, and for any $t\in (s,T],$ the distribution $ P_{s,t}^{i,x} $ of $X_{s,t}^{i,x}$ has positive density function $p_{s,t}^{i,x}$ with respect to the Lebesgue measure. When $s=0$, we simply denote
 $$X_{0,t}^{i,x}=X_t^{i,x},\ \ \ P_{0,t}^{i,x}= P_t^{i,x}.$$
 The associated Markov semigroup $(P_{s,t}^{(i)})_{0\le s\le t\le T}$  is given by
 $$P_{s,t}^{(i)} f(x):= \E[f(X_{s,t}^{i,x})],\ \ 0\le s\le t\le T, x\in \R^d, f\in \B_b(\R^d),$$
 where $\B_b(\R^d)$ is the space of all bounded measurable functions on $\R^d$.
If the initial value is random with distributions $\nu\in \scr P,$ where $\scr P$ is the set of all probability measures on $\R^d$, we denote the diffusion process  by $X_t^{i,\nu}$, which has distribution
$$P_t^{i,\nu} =\int_{\R^d} P_t^{i,x} \nu (\d x),\ \ \ i=1,2, \ t\in (0,T].$$
Let $p_t^{i,\nu}$ be the density function of $P_t^{i,\nu}$ with respect to the Lebesgue measure.

 We  intend to estimate the relative entropy
 $$\Ent(P_t^{1,\nu_1}  |P_t^{2,\nu_2}):=\int_{\R^d} \Big(\log \ff{\d P_t^{1,\nu_1} }{\d P_t^{2,\nu_2}}\Big)\d P_t^{1,\nu_1}= \E\bigg[ \Big(\log \ff{p_t^{1,\nu_1} }{p_t^{2,\nu_2}}\Big)(X_t^{1,\nu_1})\bigg]$$ for $
 t\in (0,T]$ and $\nu_1,\nu_2\in \scr P_2.$ 
Before moving on, let us recall a nice entropy inequality derived by Bogachev, R\"ockner and  Shaposhnikov \cite{BRS}.
 For a $d\times d$-matrix valued function $a=(a^{kl})_{1\le k,l\le d}$, the divergence is
 an $\R^d$-valued function defined by
  $$ {\rm div} a  := \Big(\sum_{l=1}^d \pp_l a^{kl}\Big)_{1\le k\le d},$$ where $\pp_l:=\ff{\pp}{\pp x^l}$
  for $x=(x^l)_{1\le l\le d}\in\R^d$. Let
 \beg{align*} \Phi^\nu(s,y):= &\,(a_1(s,y)-a_2(s,y)) \nn\log p_{s}^{1,\nu}(y) +{\rm div}\{a_1(s,\cdot)-a_2(s,\cdot)\}(y) \\
&  + b_2(s,y)-b_1(s,y),\ \ s\in (0,T], y\in\R^d, \nu\in \scr P,\end{align*}
where $\nn$ is the gradient operator for weakly differentiable functions on $\R^d$. In particular,
$\|\nn f\|_\infty$ is the Lipschitz constant of   $f$.

By \cite[Theorem 1.1]{BRS},  the entropy inequality
\beq\label{BOG}  \Ent(P_t^{1,\nu}|P_t^{2,\nu})\le \ff 1 2 \int_0^t    \E\big[\big|a_2(s,X_s^{1,\nu})^{-\ff 1 2 } \Phi^\nu(s,X_s^{1,\nu})\big|^2\big]\d s,\ \ t\in (0,T]\end{equation}
holds under the following   assumption $(H)$.

\beg{enumerate}\item[$(H)$]   For each $i=1,2$, $b_i$ is locally bounded, and there exists a constant $K>1$ such that
$$ \|a_i(t,x)\|\lor \|a_i(t,x)^{-1}\|\lor \|\nn a_i(t,\cdot)(x)\|\le K,\ \ (t,x)\in [0,T]\times \R^d.$$
Moreover, at least one of the following conditions hold:

$(1)$ $\int_{0}^T \E\big[\ff{\|a_2(t,X_t^{1,\nu})\|}{1+|X_t^{1,\nu}|^2}+ \ff{|b_2(t,X_t^{1,\nu})|+|\Phi^\nu(t,X_t^{1,\nu})|}{1+|X_t^{1,\nu}|}\big]\d t <\infty;$

$(2)$ there exist  $1\le V\in C^2(\R^d)$ with $V(x)\to\infty$ as $|x|\to\infty$, and a constant $K>0$ such that
$$L_{t}^{a_2,b_2}V(x)\le KV(x),\ \ \int_0^T\E\Big[\ff{|\<\Phi^\nu(t,X_t^{1,\nu}), \nn V(X_t^{1,\nu})\>|}{V(X_t^{1,\nu})}\Big]\d t<\infty.$$
\end{enumerate} It is well known that $(H)$ implies the existence and uniqueness of the diffusion processes $(X_t^{i,\nu})_{i=1,2}$ for any $\nu\in \scr P,$  and the existence of the density functions $(p_t^{i,\nu})_{i=1,2},$ see for instance \cite{BKRS}.

As observed in \cite[Remark 1.4]{BRS} that
one may  have
$$\int_0^t\E\big[  |\nn\log p_{s}^{1,\nu}|^2(X_s^{1,\nu})\big]\d s<\infty,$$
provided   $\nu$ has finite information entropy, i.e.  $\rr(x):= \ff{\d\nu}{\d x}$  satisfies $\int_{\R^d} (\rr|\log \rr|)(x)\d x<\infty.$
In this case,  \eqref{BOG} provides a non-trivial upper bound for
$\Ent(P_t^{1,\nu}|P_t^{2,\nu}).$

However, when $X_s^{1,x}$ is the standard Brownian motion starting from   a fixed initial value $x$, i.e. $\nu=\dd_x$, we have
$$\E[ |\nn\log p_{s}^{1,x}|^2(X_s^{1,x})]= \ff 1 {s^2} \E |X_s^{1,x}-x|^2 = \ff 1 s.$$
So, for elliptic diffusions the best possible short time estimate
on $\E[ |\nn\log p_s^{1,x}|^2(X_s^{1,x})]$ behaves like $ s^{-1}$, so that
$$  \int_0^t \E[ |\nn\log p_s^{1,x}|^2(X_s^{1,x})]\d s=\infty,\ \ t>0.$$
Consequently, the estimate \eqref{BOG} becomes trivial   when
\beq\label{AS2}\inf_{(s,x)\in [0,T]\times \R^d} \|a_1(s,x)-a_2(s,x)\|>0.\end{equation}

So, the key point of the present study is to cancel the small time singularity in  \eqref{BOG}, which stimulates us to develop a new coupling method, i.e. the bi-coupling method in Section 2 below.

\subsection{Main idea and structure of the paper}
To kill the singularity in \eqref{BOG} for small $t>0$, in Section 2  we introduce a new technique by constructing an interpolation diffusion process which is coupled with each of the given two diffusion processes respectively, so we call it the bi-coupling argument. In Section 3 we apply the bi-coupling to estimate the entropy between two diffusion processes, and as an application, in Section 4 we establish the entropy-cost inequality \eqref{EC} for the McKean-Vlasov SDE \eqref{E1}.

To measure the singularity/regularity of coefficients in \eqref{E1}, we introduce the following class
of   Dini functions
 $$\D:=\bigg\{\varphi: [0,\infty)\to [0,\infty)\ \text{is\ increasing\ and \ concave,}\ \varphi(0)=0,
 \int_0^1\ff{\varphi(s)}s\d s<\infty\bigg\}.$$
For $\varphi \in\D,t>0$ and a function $f$ on $[0,t]\times \R^d$, let
\beg{align*} &\|f\|_{t,\infty}:=\sup_{x\in\R^d}|f(t,x)|,\ \ \|f\|_{r\to t,\infty}:=\sup_{s\in [r,t]}\|f\|_{s,\infty},\ \ r\in [0,t],\\
&\|f\|_{0\to T,\varphi}:= \sup_{t\in [0,T],x\ne y\in \R^d} \bigg(|f(t,x)|+\ff{|f(t,x)-f(t,y)|}{\varphi(|x-y|)}\bigg).\end{align*}

In the following, $c=c(K,T,d,\varphi)$ stands for   a constant depending only on $K,T,d$ and $\varphi$ given in  $(A_1)$ and $(A_2)$.

\section{Bi-coupling method and   density estimates}

Let $\si_i=\ss{2a_i}, i=1,2.$ Consider SDEs:
 \beq\label{SDE} \d X_t^i= b_i(t,X^i_t)\d t+\si_i(t,X_t^i)\d W_t,\ \ t\in [0,T],\ i=1,2.\end{equation}  We make the following assumptions $(A_1)$ and $(A_2)$ where  $b_i$ may have a Dini continuous term with respect to some $\varphi\in \D$.

\beg{enumerate} \item[$(A_1)$] For each $i=1,2$, $b_i=b_i^{(0)}+b_i^{(1)}$ is locally bounded, and  there exists   a constant  $K>0$ such that
 \beg{align*} \|b_i^{(0)}\|_{0\to T,\infty}\lor\|\nn b_i^{(1)}\|_{0\to T,\infty} \lor\|a_i \|_{0\to T,\infty}\lor\|a_i^{-1}\|_{0\to T,\infty}\lor\|\nn a_i\|_{0\to T,\infty}\le K.\end{align*}
  \item[$(A_2)$] There exist  $i\in \{1,2\}$ and $\varphi\in \D$ such that
 $ \|b_i^{(0)}\|_{0\to T,\varphi} \le K.$
 \end{enumerate}

According to \cite[Theorem 2.1]{Ren},   $(A_1)$ implies the well-posedness of \eqref{SDE}.
For any $s\in [0,T)$ and $x\in \R^d$, let $X_{s,t}^{i,x}$ be the unique solution of \eqref{SDE}  for $t\in [s,T]$ with $X_{s,s}^{i,x}=x$. Then
  $(X_{s,t}^{i,x})_{(t,x)\in [0,T]\times\R^d}$ is the   diffusion process generated by $(L_t^{a_i,b_i})_{t\in [s,T]},\ i=1,2.$

  For fixed   $x_1,x_2\in\R^d$, let $X_t^{i,x_i}:=X_{0,t}^{i,x_i}$ solve \eqref{SDE} for $X_0^{i,x_i}=x_i.$ We have
$$P_t^{i,x_i}:=\L_{X_t^{i,x_i}},\ \ i=1,2,\ t\in (0,T].$$
To estimate $\Ent(P_{t_1}^{1,x_1}|P_{t_1}^{2,x_2})$ for some $t_1\in (0,T]$,   we choose $t_0\in (0,\ff 1 2 t_1]$ and construct a bridge diffusion process $X_t^{\<t_0\>x_1}$
starting at $x_1$ which is generated by $L_t^{a_1,b_1}$ for $t\in [0,t_0]$ and $L_t^{a_2,b_2}$ for $t\in (t_0,t_1].$ More precisely, let
\beg{align*}&b^{\<t_0\>}(t,\cdot):= 1_{[0,t_0]}(t) b_1(t,\cdot)+ 1_{(t_0,t_1]}(t)b_2(t,\cdot),\\
 &\si^{\<t_0\>}(t,\cdot):= 1_{[0,t_0]}(t) \si_1(t,\cdot)+ 1_{(t_0,t_1]}(t)\si_2(t,\cdot),\ \ \ t\in [0,t_1].\end{align*}
We consider the interpolation SDE
\beq\label{SDE'} \d X_t^{\<t_0\>x_1} = b^{\<t_0\>}(t,X_t^{\<t_0\>x_1})\d t+\si^{\<t_0\>}(t,X_t^{\<t_0\>x_1})\d W_t,\ \ \ X_0^{x_1}=x_1,\ t\in [0,t_1].\end{equation}
Let $P_t^{\<t_0\>x_1}:=\L_{X_t^{\<t_0\>x_1}}.$
We will deduce  from    \eqref{BOG} a finite upper bound for $\Ent(P_{t_1}^{1,x_1}|P_{t_1}^{\<t_0\>x_1}),$
 where the singularity at $t=0$ disappears since the distance of diffusion coefficients vanishes for $t\in [0,t_0].$ Moreover, we will estimate the moment for the density of  $  P_{t_1}^{\<t_0\>x_1}$ with respect to $P_{t_1}^{2,x_2},$ so that by the following entropy inequality \eqref{EP},  we derive the desired upper bound on $\Ent(P_{t_1}^{1,x_1}|P_{t_1}^{2,x_2})$.  We remark that \eqref{EP} has been presented in \cite{HS} for $p=2,$
but in the present study  we  shall need the inequality for  $p>2$ as required in the dimension-free Harnack inequality  due to \cite{Ren}, see the proof of Proposition \ref{P2.1} for details.

\beg{lem}\label{LN0} Let $\mu_1,\mu_2$ and $\mu$ be probability measures on a measurable space $(E,\B).$
Then for any $p>1$,
\beq\label{EP} \Ent(\mu_1|\mu_2)\le p\Ent(\mu_1|\mu)+  (p-1)\log \int_E \Big(\ff{\d \mu}{\d\mu_2}\Big)^{\ff{p}{p-1}}\d\mu_2,\end{equation}
where the right hand side is set to be infinite if $\ff{\d\mu_1}{\d\mu}$ or $\ff{\d\mu}{\d\mu_2}$
does not exist. \end{lem}
\beg{proof} It suffices to prove for the case that $\ff{\d\mu_1}{\d\mu}$ and $\ff{\d\mu}{\d\mu_2}$ exist
such that the upper bound is finite.
In this case, we have
\beg{align*} &\Ent(\mu_1|\mu_2)-\Ent(\mu_1|\mu)= \int_E\Big\{\log \ff{\d\mu_1}{\d\mu_2}- \log
\ff{\d\mu_1}{\d\mu}\Big\}\d\mu_1 \\
 &= \int_E \Big\{\log\ff{\d\mu}{\d\mu_2}\Big\}\d\mu_1
  = \ff {p-1}p \int_E  \Big(\ff{\d\mu_1}{\d\mu_2}\Big)\log\Big(\ff{\d\mu}{\d\mu_2}\Big)^{\ff p{p-1}}\d\mu_2.\end{align*}
 Combining  with the Young inequality  \cite[Lemma 2.4]{SPA09}, we obtain
 $$\Ent(\mu_1|\mu_2)-\Ent(\mu_1|\mu)\le \ff {p-1}p \Ent(\mu_1|\mu_2)
 +\ff {p-1}p \log \int_E \Big(\ff{\d\mu}{\d\mu_2}\Big)^{\ff p{p-1}} \d\mu_2.$$
  \end{proof}

By Lemma \ref{LN0}, for any $p>1$ we have
\beq\label{YPP} \Ent(P_{t_1}^{1,x_1}|P_{t_1}^{2,x_2})\le p \Ent(P_{t_1}^{1,x_1}|P_{t_1}^{\<t_0\>x_1})+(p-1)\log \int_{\R^d} \bigg(\ff{\d P_{t_1}^{\<t_0\>x_1}}{\d P_{t_1}^{2,x_2}}\bigg)^{\ff{p}{p-1}}\d P_{t_1}^{2,x_2}.\end{equation}
Noting that $a(t,\cdot)- a_1(t,\cdot)=0$ for $t\in [0,t_0],$  we may apply \eqref{BOG} to derive a non-trivial upper bound on the first term in the right hand side of \eqref{YPP}, see Proposition \ref{P3.1} for details. So, in the following, we only estimate the second term. To this end, we need the following simple lemma.

\beg{lem}\label{LLN} Let $\xi_t\ge 0$ be a continuous semi-martingale such that
$$\d \xi_t\le k_1\xi_t\d t+\d A_t+\d M_t,\ \ t\in [0,T],$$
where $k_1>0$ is a constant, $A_t$ is an  increasing function with $A_0=0$, and $M_t$ is a local martingale with
$$\d\<M\>_t\le k_1\xi_t\d t.$$
Then for any $t_0\in (0, T\land k_1^{-1})$ and constants $\ll,k>0$ such that
\beq\label{R*} k (1-k_1t_0)\ge k_1\Big(1+\ff \ll 2 \Big),\end{equation}
we have
$$\E\exp\Big[\ff{\ll\xi_{t_0}}{1+kt_0}\Big]\le \exp\big[\ll\xi_0+\ll A_{t_0}\big].$$
\end{lem}

\beg{proof} Let $\eta_t:= \exp\big[\ff{\ll\xi_{t}}{1+kt}\big].$ By It\^o's formula, we find a local martingale $\tt M_t$ such that
\beg{align*} &\d \eta_t=\eta_t\Big\{\ff{\ll}{1+kt}\d\xi_t+\ff{\ll^2}{2(1+kt)^2}\d\<M\>_t -\ff{k\ll \xi_t}{(1+kt)^2} \d t\Big\}+\d\tt M_t\\
&\le \eta_t\xi_t \Big\{\ff{\ll k_1}{1+kt} +\ff{\ll^2k_1}{2(1+kt)^2}  -\ff{k\ll  }{(1+kt)^2} \Big\}\d t+\ll \eta_t \d A_t+\d\tt M_t,\ \ \ t\in [0,T].\end{align*}
By \eqref{R*} we have
$$\ff{\ll k_1}{1+kt} +\ff{\ll^2k_1}{2(1+kt)^2}  -\ff{k\ll  }{(1+kt)^2} \le 0,\ \ t\in [0,t_0],$$
so that
$$\d \eta_t\le \ll \eta_t\d A_t+\d \tt M_t,\ \ t\in [0,t_0].$$
By Gronwall's lemma, this implies
$$\E[\eta_{t_0}] \le \eta_0\e^{\ll A_{t_0}},$$
which coincides with the desired estimate.

\end{proof}

\beg{prp}\label{P2.1} Assume $(A_1)$ and $(A_2)$. Then there exist  constants $p=p(K,T,d)>2, \vv =\vv(K,T,d)\in (0,\ff 1 2]$ and $c=c(K,T,d)>0,$ such that for any $x_1,x_2\in\R^d, t_1\in (0,T]$ and $t_0=\vv t_1,$
$$\log  \int_{\R^d} \Big(\ff{\d P_{t_1}^{\<t_0\>x_1}}{\d P_{t_1}^{2,x_2}}\Big)^{\ff{p}{p-1}}\d P_{t_1}^{2,x_2}
\le \ff c{t_1}\bigg(|x_1-x_2|^2  + \int_{0}^{t_1}\big\{\|a_1-a_2\|_{t,\infty}^2 + \|b_1-b_2\|_{t,\infty}^2\big\}\d t\bigg).$$
\end{prp}

\beg{proof}
(a) Recall that $\B_b(\R^d)$ is the space of all bounded measurable functions on $\R^d$, and let
$$P_t^{\<t_0\>} f(x):=\E[f(X_t^{\<t_0\>x})],\ \ P_t^{(2)}f(x):= \E[f(X_t^{2,x})],\ \ f\in \B_b(\R^d),\ (t,x)\in [0,T]\times \R^d.$$
Then  the desired estimate follows from the inequality 
\beq\label{WT} \beg{split}&\big|P_{t_1}^{\<t_0\>} f(x_1)\big|^p\le \big(P_{t_1}^{(2)}|f|^p(x_2)\big)\\
&\times \exp\bigg[\ff{c(p-1)}{t_1}\bigg(|x_1-x_2|^2  + \int_{0}^{t_1}\big\{\|a_1-a_2\|_{t,\infty}^2 + \|b_1-b_2\|_{t,\infty}^2\big\}\d t\bigg)\bigg],\ \ f\in \B_b(\R^d).\end{split}\end{equation}
Indeed,   taking $f:= \big(n\land  \ff{\d P_{t_1}^{\<t_0\>x_1}}{\d P_{t_1}^{2,x_2}}\big)^{\ff{1}{p-1}}$ for $n\ge 1$, this inequality implies
\beg{align*} &\bigg(\int_{\R^d} \Big(n\land  \ff{\d P_{t_1}^{\<t_0\>x_1}}{\d P_{t_1}^{2,x_2}}\Big)^{\ff{p}{p-1}} \d P_{t_1}^{2,x_2}\bigg)^p
\le \bigg(\int_{\R^d} \Big(n\land  \ff{\d P_{t_1}^{\<t_0\>x_1}}{\d P_{t_1}^{2,x_2}}\Big)^{\ff{1}{p-1}} \d P_{t_1}^{\<t_0\>x_1}\bigg)^p\\
&\le \bigg(\int_{\R^d} \Big(n\land  \ff{\d P_{t_1}^{\<t_0\>x_1}}{\d P_{t_1}^{2,x_2}}\Big)^{\ff{p}{p-1}} \d P_{t_1}^{2,x_2}\bigg)  \\
&\quad \times \exp\bigg[\ff{c(p-1)}{t_1}\bigg(|x_1-x_2|^2  + \int_{0}^{t_1}\big\{\|a_1-a_2\|_{t,\infty}^2 + \|b_1-b_2\|_{t,\infty}^2\big\}\d t\bigg)\bigg].\end{align*}
Taking $\log$ in both sides we derive 
\beg{align*} & \log  \int_{\R^d} \Big(n\land \ff{\d P_{t_1}^{\<t_0\>x_1}}{\d P_{t_1}^{2,x_2}}\Big)^{\ff{p}{p-1}}\d P_{t_1}^{2,x_2}\\
&\le \ff c{t_1}\bigg(|x_1-x_2|^2  + \int_{0}^{t_1}\big\{\|a_1-a_2\|_{t,\infty}^2 + \|b_1-b_2\|_{t,\infty}^2\big\}\d t\bigg),\ \ n\ge 1,\end{align*} 
which implies the desired estimate as
    $n\to\infty$. So, it remains to find constants $p>2$ and $c>0$ such that \eqref{WT} holds. 
 
Let $(P_{s,t}^{(2)})_{0\le s\le t\le T}$ be the semigroup generated by $L_t^{a_2,b_2}$, i.e.
$$P_{s,t}^{(2)}f(x):= \E[f(X_{s,t}^{2,x})],\ \ \ f\in \B_b(\R^d),$$
where $(X_{s,t}^{2,x})_{t\in [s,T]}$ solves
$$\d X_{s,t}^{2,x}= b_2(t, X_{s,t}^{2,x})\d t + \si_2(t,X_{s,t}^{2,x})\d W_t,\ \ X_{s,s}^{2,x}=x,\ t\in [s,T].$$
By the Markov property and the SDE \eqref{SDE'}, we obtain
\beq\label{SM} P_{t_1}^{\<t_0\>}f(x_1)= \E \big[(P_{t_0,t_1}^{(2)} f)(X_{t_0}^{1,x_1})\big],\ \
P_{t_1}^{(2)}f(x_2)= \E \big[(P_{t_0,t_1}^{(2)} f)(X_{t_0}^{2,x_2})\big].\end{equation}
 By \cite[Theorem 2.2]{Ren} which applies to a more general setting where $b_2^{(0)}$ only satisfies a local integrability condition,
 and noting that $t_1-t_0= (1-\vv)t_1$, 
we find constants $p_1=p_1(K,T,d)> 1\lor \ff d 2$ and $c_1=c_1(K,T,d,\vv)>0$ such that
 \beq\label{SM2} \big|P_{t_0, t_1}^{(2)} f(x)\big|^{p_1}\le \big(P_{t_0,t_1}^{(2)}|f|^{p_1} (y)\big)
 \e^{\ff{c_1|x-y|^2}{t_1}},\ \ f\in\B_b(\R^d), x,y\in\R^d.\end{equation}
Combining this with \eqref{SM} and Jensen's inequality, for $p:=2p_1>2\lor d$  we obtain
\beq\label{SM3} \beg{split}&\big|P_{t_1}^{\<t_0\>}f(x_1)|^{p} = \big|\E[P_{t_0,t_1}^{(2)}f(X_{t_0}^{1,x_1})]\big|^{2p_1}
\le \Big(\E\big[|P_{t_0,t_1}^{(2)}f|^{p_1}(X_{t_0}^{1,x_1})\big]\Big)^2\\
& \le \bigg\{\E\Big[\big(P_{t_0,t_1}^{(2)}|f|^{p_1}(X_{t_0}^{2,x_2})\big) \exp\Big(\ff{c_1|X_{t_0}^{1,x_1}-X_{t_0}^{2,x_2}|^2}{t_1}\Big)\Big]\bigg\}^2\\
&\le \big(\E\big[P_{t_0,t_1}^{(2)} |f|^{2p_1}(X_{t_0}^{2,x_2})\big] \big)\E\bigg[\exp\Big(\ff{2c_1|X_{t_0}^{1,x_1}-X_{t_0}^{2,x_2}|^2}{t_1}\Big)\bigg]\\
& = \big(P_{t_1}^{(2)}|f|^p(x_2)\big) \E\bigg[\exp\Big(\ff{2c_1|X_{t_0}^{1,x_1}-X_{t_0}^{2,x_2}|^2}{t_1}\Big)\bigg].
\end{split}\end{equation}
Thus, to prove \eqref{WT},  it remains to estimate the  expectation term in the upper bound.

 (b) Since   the exponential term is symmetric in $(X_{t_0}^{1,x_1},X_{t_0}^{2,x_2}),$ without loss of generality, in $(A_2)$ we may and do assume that $\|b_1^{(0)}\|_{0\to T,\varphi}\le K.$
  We shall use Zvonkin's transform to kill this non-Lipschitz term.
   By $(A_1)$, $b_1^{(0)}$ is bounded, and noting that $p:=2p_1>2\lor d$, for a fixed  constant $q>2$ such that $\ff d p+\ff 2 q<1$, we have $\|b_1^{(0)}\|_{\tt L_q^p}<\infty.$    
 So, according to   \cite[Theorem 2.1]{ZY},   there exist constants $c_1=c_1(K,T,d,p,q)>0$ and
$\bb=\bb(p,q)\in (0,1)$ such that for any $\ll>0$, the PDE
\beq\label{W3} (\pp_t +L_t^{a_1,b_1}-\ll)u_t= -b_1^{(0)}(t,\cdot),\ \ \ t\in [0,T], u_T=0\end{equation}
has a unique
solution satisfying
\beq\label{FM0}\ll^\bb (\|u\|_{0\to T,\infty}+\|\nn u\|_{0\to T,\infty})+ \|\pp_t u\|_{\tt L_q^p}+\|\nn^2 u\|_{\tt L_q^p}\le c_1,\end{equation}
where for any measurable function $g$ on $[0,T]\times \R^d$, 
\beq\label{LQP} \|g\|_{\tt L_q^p}:= \sup_{z\in\R^d} \bigg(\int_0^T\|1_{B(z,1)} g(t,\cdot)\|_{L^p(\R^d)}^q\d t\bigg)^{\ff 1 q}.\end{equation}  Let $P_{s,t}^{a_1,b_1^{(1)}}$ be the Markov semigroup generated by $L_t^{a_1,b_1^{(1)}}$, and let
$p_{s,t}^{a_1,b_1^{(1)}}$ be the heat kernel with respect to the Lebesgue measure.
By Duhamel's formula, we have
\beq\label{FM} u_s=\int_s^T \e^{-\ll (t-s)}P_{s,t}^{a_1,b_1^{(1)}} \big\{\nn_{b_1^{(0)}} u_t+b_1^{(0)}(t,\cdot)\big\}\d t,\ \ s\in [0,T].\end{equation}
Let $\nn_x^2$ be the Hessian operator in $x$. By \cite[Theorem 1.2]{MPZ}, under $(A_1)$ we find a constant $\dd=\dd(K,T,d)>1$ such that
$$|\nn^2_x p_{s,t}^{a_1,b_1^{(1)}}(x,y)|\le \ff {\ll}{t-s} g_\dd(t-s, x,y),\ \ 0\le s<t\le T, x,y\in \R^d$$
holds for
\beq\label{GD} g_\dd(r,x,y):= (\pi\dd r)^{-\ff d 2}\e^{-\ff{|\theta_{s,t}(x)-y|^2}{\dd r}},\ \ r>0, x,y\in\R^d,\end{equation}
where $\theta: [0,T]\times [0,T]\times \R^d\to\R^d$ is a measurable map.
So,  letting
\beq\label{HPT} h_{t}(y):=  (b_1^{(0)}(t,y)\cdot\nn) u_t(y) + b_1^{(0)}(t,y),\end{equation}
and denoting by $(\nn_x,\nn_x^2)$   the gradient and Hessian operators in $x\in\R^d$, we obtain
\beq\label{HPT2}\beg{split} &|\nn^2_x u_s(x)|\le \int_s^T\ff{ \e^{-\ll(t-s)}}{t-s}\big|\nn^2_x  P_{s,t}^{a_1,b_1^{(1)}}(h_t-h_t(z))(x)\big|_{z=\theta_{s,t}(x)}\d t\\
&\le \int_s^T\ff{\e^{-\ll(t-s)}}{t-s}\d t\int_{\R^d}  \big|\nn^2_x p_{s,t}^{a_1,b_1^{(1)}}(x,y)|\cdot |h_t(y)-h_t(\theta_{s,t}(x))|\d y.\end{split}\end{equation}
By $(A_2)$, \eqref{FM0} for $\ll\ge 1$, and \eqref{HPT}, we have
\beq\label{HPT3} |h_t(y)-h_t(\theta_{s,t}(x))|\le (1+c_1)|b_1^{(0)}(t,y)- b_1^{(0)}(t,\theta_{s,t}(x))|+ K|\nn u_t(y)-\nn u_t(\theta_{s,t}(x))|.\end{equation}
In the following, we estimate these two terms in the upper bound respectively.

Since $\varphi$ is concave, we find a constant $c_2=c_2(K,T,d)>0$ such that
 \beg{align*} & \int_{\R^d} |b_1^{(0)}(t,y)-b_1^{(0)}(t,\theta_{s,t}(x))|g_\dd(t-s, x,y)\d y\\
&\le
K\int_{\R^d} \varphi(|y-\theta_{s,t}(x)|)  g_\dd(t-s, x,y)\d y\\
&\le
K\varphi\bigg(\int_{\R^d}  |y-\theta_{s,t}(x)|  g_\dd(t-s, x,y)\d y\bigg)
\le c_2 \varphi\Big(\ss{t-s}\Big),\ \ 0\le s<t\le T,x\in \R^d.\end{align*}
Hence,
\beq\label{B1} \beg{split}&\sup_{s\in [0,T] } \int_s^T\ff{\e^{-\ll (t-s)}}{t-s}\d t \int_{\R^d} |b_1^{(0)}(t,y)-b_1^{(0)}(t,\theta_{s,t}(x))|g_\dd(t-s, x,y)\d y\\
&\le c_2 \int_0^T \ff{\e^{-\ll t }\varphi(t^{\ff 1 2})}{t} \d t=:  \vv_1,\end{split}\end{equation}
where  $\vv_1=\vv_1(\ll, K,T,d,\varphi). $  Since $\varphi\in\D$ implies 
$$\int_0^T \ff{\varphi(t^{\ff 1 2})}{t} \d t =2\int_0^{T^{\ff 1 2} } \ff{\varphi(s)}s \d s<\infty,$$ 
by the dominated convergence theorem we derive
$\lim_{\ll\to\infty} \vv_1=0.$

On the other hand, let $\aa=1-\ff d p\in (0,1)$. By the Sobolev embedding theorem, see e.g. \cite{AF},
there exists a constant   $c_0>0$ depending on $p$ and $d$ such that
$$ \sup_{z\ne y\in B(z,1)}\ff{|f(y)-f(z)|}{|y-z|^\aa}\le c_0\|1_{B(z,1)}(|f|+|\nn f|)\|_{L^p},\ \ z\in\R^d,\ \ f\in W_{loc}^{1,p}(\R^d).$$
  So,
  $$|\nn u_t(y)-\nn u_t(z)|\le c_0|y-z|^\aa \|1_{B(z,1)}(|\nn u_t|+\|\nn^2 u_t\|)\|_{L^p(\R^d)}\big),\ \ \text{if}\ |y-z|<1.$$
 Noting that $\ff d p+\ff 2 q<1$ and $\aa= 1-\ff d p$ imply $(1-\aa)\ff q {q-1}<1$,  by
 combining this with \eqref{FM0} and \eqref{GD}, we find   constants $c_3=c_3(p,d)>0$ and   $\vv_2=\vv_2(\ll,K,T,d,p,q)>0$, where  $\vv_2\to 0$ as $\ll\to\infty$, such that
 \beg{align*}& \int_s^T \ff{ \e^{-\ll(t-s)}}{t-s}\d t\int_{\R^d} |\nn u_t(y)-\nn u_t(\theta_{s,t}(x))|g_\dd(t-s, x,y)\d y\\
 &\le c_3 \bigg(\int_s^T \e^{-\ll (t-s)} (t-s)^{-(1-\aa)\ff q{q-1}}\d t \bigg)^{\ff {q-1}q}\big(\|\nn u\|_{0\to T,\infty}+\|\nn^2u\|_{\tt L_q^p}\big) \le \vv_2,\ \ s\in [0,T].\end{align*}
By \eqref{FM0}, and combining this with  \eqref{HPT2}, \eqref{HPT3}, and \eqref{B1}, we find large enough $\ll=\ll(K,T,P,\varphi)>0$ such that $\|\nn^2 u\|_{0\to T,\infty}\le \ff 1 2.$ Combining this with
\eqref{FM0}, we may choose large enough $\ll>0$ such that
\beq\label{W4} \|u\|_{0\to T,\infty}\lor \|\nn u\|_{0\to T,\infty}\lor \|\nn^2 u\|_{0\to T,\infty}\le \ff 1 {2}.\end{equation}
In particular, letting
\beq\label{TT} \tt X_t^{i,x_i}:= X_t^{i,x_i}+u_t(X_t^{i,x_i}),\ \ \ i=1,2,\end{equation}  we have
 \beq\label{W5} \ff 1 2|X_t^{1,x_1}-X_t^{2,x_2}|\le |\tt X_t^{1,x_1}-\tt X_t^{2,x_2}|\le 2 |X_t^{1,x_1}-X_t^{2,x_2}|.\end{equation}
 Hence, to bound the exponential moment in \eqref{SM3}, it suffices to estimate the corresponding term
 for $|\tt X_{t_0}^{1,x_1}-\tt X_{t_0}^{2,x_2}|^2$ replacing $|X_{t_0}^{1,x_1}-X_{t_0}^{2,x_2}|^2.$

 (c) Let $I_d$ be the $d\times d$ identity matrix. By \eqref{W3}, \eqref{TT} and It\^o's formula, we obtain
\beq\label{SV1} \beg{split}&\d \tt X_t^{1,x_1}= \big\{\ll u_t+b_1^{(1)}(t,\cdot)\big\}(X_t^{1,x_1})\d t
 + \big\{I_d +\nn u_t(X_t^{1,x_1})\big\}\si_1(t,X_t^{1,x_1})\d W_t,\\
 &\d \tt X_t^{2,x_2}= \big\{\ll u_t+ (L_t^{a_2,b_2}-L_t^{a_1,b_1})u_t + (b_2-b_1^{(0)})(t,\cdot)\big\}
(X_t^{2,x_2})\d t\\
&\qquad\qquad\qquad\qquad
 + \big\{I_d +\nn u_t(X_t^{2,x_2})\big\}\si_2(t,X_t^{2,x_2})\d W_t.\end{split}\end{equation}
 By $(A_1)$, \eqref{W4}, \eqref{W5}, and It\^o's formula, we find $k_1=k_1(K,T,d,\varphi)>0$ such that
 \beq\label{W6} \d |\tt X_{t}^{1,x_1}-\tt X_{t}^{2,x_2}|^2\le k_1 |\tt X_{t}^{1,x_1}-\tt X_{t}^{2,x_2}|^2  \d t+\d A_t + \d M_t,\ \ t\in [0,t_0],\end{equation}
where
\beq\label{AT} A_t:= k_1\int_0^t (\|a_1-a_2\|_{s,\infty}^2+\|b_1-b_2\|_{s,\infty}^2\big)\d s,\end{equation}  and
$M_t$ is a martingale satisfying
\beq\label{W7} \d\<M\>_t \le k_1 |\tt X_{t}^{1,x_1}-\tt X_{t}^{2,x_2}|^2 \d t.\end{equation}
For any $n\ge 1,$ let
$$\tau_n:=t_0\land \inf\big\{t\ge 0:\ |\tt X_{t}^{1,x_1}-\tt X_{t}^{2,x_2}|\ge n\big\},
\ \ \gg_n:= \sup_{t\in [0,\tau_n]} |\tt X_{t}^{1,x_1}-\tt X_{t}^{2,x_2}|^2.$$
By \eqref{W5} we have
\beq\label{*NN}|\tt X_{0}^{1,x_1}-\tt X_{0}^{2,x_2}|^2\le 4 |x_1-x_2|^2.\end{equation}
Moreover, to apply Lemma \ref{LLN}, let
$$t_0:= \ff{t_1}{2[(Tk_1+4k_1c_1)\lor 1]},\ \ \ll:= \ff{8c_1(1+kt_0)}{t_1},
\ \ k=\ff{k_1}{1-k_1 t_0}\Big(1+\ff \ll 2\Big),$$
so that \eqref{R*} holds and
$$\ff \ll{1+kt_0} = \ff{8c_1}{t_1}.$$
Combining this with \eqref{W6}-\eqref{*NN}, we may apply Lemma \ref{LLN} for $\xi_t=|\tt X_t^{1,x_1}-\tt X_t^{2,x_2}|^2$  to find a constant $k_2=k_2(K,T,d,\varphi)>0$ such that
$$\E\big[\e^{\ff{8c_1}{t_1}|\tt X_{t_0}^{1,x_1}-\tt X_{t_0}^{2,x_2}|^2}\big]\le \e^{\ff{k_2}{t_1}\big(|x_1-x_2|^2+\int_0^{t_0}  (\|a_1-a_2\|_{t,\infty}^2+  \|b_1-b_2\|_{t,\infty}^2)\d t\big)}.$$
This together with \eqref{SM3} implies \eqref{WT} for some constant $c=c(K,T,d,\varphi)$, and hence finishes the proof.

\end{proof}

 \section{Entropy estimates between two diffusion processes}

With the bi-coupling method and density estimates addressed in Section 2, we are able to prove the following result on entropy upper bound estimates for diffusion processes with arbitrary initial distributions in $\scr P_2$, for which the existing estimates may be invalid  as explained in Section 1.2.

 \beg{thm}\label{T1} Assume   $(A_1)$ and $(A_2)$. Then the following assertions hold for    some constants $c=c(K,T,d,\varphi)>0$ and $\vv=\vv(K,T,d,\varphi)\in (0,\ff 1 2].$
 \beg{enumerate}\item[$(1)$]
 For any $\nu_1,\nu_2\in\scr P$ and $t\in (0,T],$
\beq\label{AE1} \beg{split}\Ent(P_t^{1,\nu_1}|P_t^{2,\nu_2})\le&\,  \ff {c\W_2(\nu_1,\nu_2)^2}t +  \ff c t \int_{0}^t \big\{\|b_1-b_2\|_{s,\infty}^2+\|a_1-a_2\|_{s,\infty}^2\big\}\d s\\
&+ c\Big[  \|a_1-a_2\|_{\vv t\to t,\infty}^2+
   \int_{\vv t}^t\|{\rm div}(a_1-a_2)\|_{s,\infty}^2\d s\Big].\end{split}\end{equation}
 \item[$(2)$]    If there exists a constant $C(K)>0$ such that
   \beq\label{XY} \|\nn^i b_1\|_{0\to T,\infty}+\|\nn^ia_1\|_{0\to T,\infty}\le C(K),\ \ i=1,2,\end{equation}
   then  for any $\nu_1,\nu_2\in\scr P$ and $t\in (0,T],$
   \beq\label{AE3} \beg{split} \Ent(P_t^{1,\nu_1}|P_t^{2,\nu_2})\le &\, \ff c t\bigg[ \W_2(\nu_1,\nu_2)^2+
    \int_{0}^t \big(\|b_1-b_2\|_{s,\infty}^2+\|a_1-a_2\|_{s,\infty}^2\big)\d s\bigg]\\
    &+\int_{\vv t}^t
    \|{\rm div}(a_1-a_2)\|_{s,\infty}^2  \d s.\end{split} \end{equation}
\end{enumerate}
\end{thm}
To prove Theorem \ref{T1}, we shall apply 
  \eqref{YPP}, where  the second term in the upper bound has been estimated in Proposition \ref{P2.1}, and the first term  will be estimated by using  \eqref{BOG} and   the following result.

\beg{prp}\label{P3.1} Assume $(A_1)$. Then the following assertions hold.
\beg{enumerate} \item[$(1)$] There exists a constant $c=c(K,T,d)>0$ such that
\beq\label{YD1} \int_r^t\d s\int_{\R^d}\ff{|\nn p_s^{1,x}|^2} {p_s^{1,x}} (y)\d y\le c \log \Big(1+\ff t r\Big),\ \ 0<r\le t\le T, x\in\R^d.\end{equation}
\item[$(2)$]  If   $\eqref{XY}$ holds,  then exists a constant $c=c(K,T,d)>0$ such that
\beq\label{YD3}   \int_{\R^d}\ff{|\nn p_t^{1,x}|^2} {p_t^{1,x}} (y)\d y\le \ff c  t, \  \  \ t\in (0, T],x\in\R^d.\end{equation}\end{enumerate}
\end{prp}

To prove  Proposition \ref{P3.1},
we first present the following lemma.

\beg{lem}\label{LN1} Assume $(A_1)$ with  the condition on $\|\nn a_1\|_{0\to T,\infty}$ replacing by the weaker one: there exists $\bb\in (0,1)$ such that
$$\|a_1(t,x)-a_1(t,y)\|\le K|x-y|^\bb,\ \ t\in [0,T], x,y\in \R^d.$$
Then  there exists  a  constant  $c =c (K,T,d,\bb) >0$ such that
\beq\label{L2}\beg{split}& \bigg|\int_{\R^d} (p_r^{1,x}\log p_r^{1,x})(y)\d y - \int_{\R^d} (p_t^{1,x}\log p_t^{1,x})(y)\d y\bigg|\\
&\qquad \le c \log \Big(1+\ff t r\Big),\ \ \ \ 0<r\le t\le T,
x\in \R^d.\end{split}\end{equation}

\end{lem}
\beg{proof} Let $x\in\R^d$ be fixed. Simply denote $\rr_t(y):=p_t^{1,x}(y), t\in (0,T], y\in\R^d.$
Let $\theta_t(x)$ solve
\beq\label{TH} \pp_t \theta_t(x)= b_1(t,\theta_t(x)),\ \ \theta_0(x)=x,\ \ t\in [0,T].\end{equation}
By \cite[Theorem 1.2]{MPZ},   there exists a
constant $c_0=c_0(K,T,d)>1$   such that
\beq\label{GPP} \ff{1}{c_0t^{\ff d 2}}  \e^{-\ff{c_0|\theta_t(x)-y|^2} {t}}\le  \rr_t(y)
 \le \ff{c_0}{t^{\ff d 2}}  \e^{-\ff{|\theta_t(x)-y|^2} {c_0t}},\ \ x,y\in\R^d, t\in (0,T]. \end{equation}
Consequently,
\beq\label{ANN} \int_{\R^d}\rr_t(y)\log \rr_t(y) \d y \le \int_{\R^d} \rr_t(y) \log[c_0t^{-\ff d 2}]\d y = \log[c_0t^{-\ff d 2}],\ \ t\in (0,T].\end{equation}
On the other hand, by \eqref{GPP} and Jensen's inequality, there exists a constant $c_1>0$ such that
\beg{align*} & -\int_{\R^d}\rr_t(y)\log \rr_t(y) \d y = 2\int_{\R^d}\rr_t(y)\log \rr_t(y)^{-\ff 1 2}  \d y \le 2\log \int_{\R^d} \rr_t(y)^{\ff 1 2}\d y \\
& \le  2\log\bigg[ c_0^{\ff 1 2}t^{-\ff d 4} \int_{\R^d} \e^{-\ff{|\theta_t(x)-y|^2} {2c_0t}}\d y\bigg] = 2 \log[c_1^{-\ff 1 2}t^{-\ff d 4}]=\log[c_1^{-1} t^{\ff d 2}].\end{align*}
Hence,
$$\int_{\R^d}\rr_t(y)\log \rr_t(y) \d y \ge \log[c_1 t^{-\ff d 2}],\ \ t\in (0,T].$$
Combining this with \eqref{ANN}, we derive find a constant $c>0$ such that
$$\int_{\R^d}\rr_r(y)\log \rr_r(y) \d y-\int_{\R^d}\rr_t(y)\log \rr_t(y) \d y \ge \log [c_1r^{-\ff d 2}] - \log [c_0t^{-\ff d 2}] \ge -c \log\Big(1+\ff t r\Big),$$
and similarly,
$$\int_{\R^d}\rr_r(y)\log \rr_r(y) \d y-\int_{\R^d}\rr_t(y)\log \rr_t(y) \d y \le \log [c_0r^{-\ff d 2}] - \log [c_1t^{-\ff d 2}] \le c \log\Big(1+\ff t r\Big).$$
So, \eqref{L2}  holds. 
\end{proof}

\beg{proof}[Proof of Proposition \ref{P3.1}]
 Let $x\in\R^d$ be fixed, and  simply denote  $\rr_t:=p_{t}^{1,x}.$

(a) We first consider the smooth case where
\beq\label{GH}\|\nn^i b_1\|_{0\to T,\infty}+\|\nn^i a_1\|_{0\to T,\infty}<\infty,\ \ i\ge 1.\end{equation}
By \cite[Theorem 1.2]{MPZ},  there exist  a constant  $\ll>1$ and a measurable map $\theta:[0,T]\to\R^d$ such that
\beq\label{B00} \ll^{-1}  t^{-\ff {d +i} 2} \e^{-\ff{\ll |\theta_t-y|^2}  t}  \le \big|\nn^i \rr_t \big|(y)\le \ll t^{-\ff {d +i}2} \e^{-\ff{|\theta_t-y|^2}{\ll t}}, \ \ t\in (0,T],   y\in\R^d, i=0,1,2.\end{equation} Moreover, by the   Kolmogorov forward equation and integration by parts formula, we have
\beq\label{B10} \pp_t\rr_t= {\rm div }\Big[a_1(t,\cdot) \nn\rr_t+\rr_t \{{\rm div }a_1(t,\cdot)-b_1(t,\cdot)\}\Big],\ \ t\in (0,T].\end{equation}
  By   \eqref{B00},   \eqref{B10} and integration by parts formula,  we obtain
\beq\label{GJ0}\beg{split}& \int_{\R^d} \big\{\rr_t\log \rr_t - \rr_r\log \rr_r \big\}(y)\d y  = \int_r^t\d s \int_{\R^d} \big\{(1+\log \rr_s)\pp_s\rr_s   \big\}(y)\d y  \\
&= -\int_r^t\d s \int_{\R^d} \Big\<a_1(s,\cdot)\nn\log\rr_s+{\rm div}a_1(s,\cdot)-b_1(s,\cdot), \nn\rr_s\Big\> (y)\d y.\end{split}\end{equation}
Since  $a_1\ge K^{-1} I_d$, this implies
\beq\label{GJ0'}\beg{split} & \int_{\R^d} \big\{\rr_t\log \rr_t - \rr_r\log \rr_r \big\}(y)\d y +  \ff 1 K \int_r^t\d s \int_{\R^d} \ff{|\nn \rr_s|^2}{  \rr_s} (y)\d y\\
&\le - \int_r^t\d s \int_{\R^d} \Big\<{\rm div}a_1(s,\cdot)-b_1(s,\cdot), \nn\rr_s \Big\>(y)\d y\\
&=\int_r^t\d s \int_{\R^d} \Big\<\big[b_1^{(0)}-{\rm div}a_1\big](s,\cdot), \nn\rr_s \Big\>(y)\d y+ \int_r^t\d s \int_{\R^d} \Big\<b_1^{(1)}(s,\cdot), \nn\rr_s \Big\>(y)\d y.\end{split}\end{equation}
By \eqref{GH},  \eqref{B00} and Lemma \ref{LN1}, we derive
\beq\label{FN} \int_r^t\d s \int_{\R^d} \ff{|\nn \rr_s|^2}{  \rr_s} (y)\d y<\infty.\end{equation}
Noting that $(A_1)$ implies $|b_1^{(0)}-{\rm div}a_1|\le 2K$, so that
\beg{align*}&\int_r^t\d s \int_{\R^d} \Big\<\big[b_1^{(0)}-{\rm div}a_1\big](s,\cdot), \nn\rr_s \Big\>(y)\d y\\
&\le \ff 1 {2K} \int_r^t\d s \int_{\R^d} \ff{|\nn\rr_s|^2}{\rr_s}(y)\d y +   2K^3\int_r^t\d s \int_{\R^d}\rr_s(y)\d y\\
&
= \ff 1 {2K} \int_r^t\d s\int_{\R^d} \ff{|\nn\rr_s|^2}{\rr_s}(y)\d y +  2 K^3(t-r).\end{align*}
Moreover, by the integration by parts formula, \eqref{B00}  and $\|\nn b_1^{(1)}\|_{0\to T,\infty}\le K$,
we obtain
$$\int_r^t\d s \int_{\R^d} \Big\<b_1^{(1)}(s,\cdot), \nn\rr_s \Big\>(y)\d y=-\int_r^t\d s \int_{\R^d}
{\rm div} \{b_1^{(1)}(s,y)\} \rr_s(y)\d y\le K(t-r).$$
Combining these with \eqref{GJ0'} and \eqref{FN}, we derive
\beq\label{ADD}\beg{split}&\int_r^t\d s \int_{\R^d} \ff{|\nn \rr_s|^2}{  \rr_s} (y)\d y\\
&\le 2K \int_{\R^d} \big\{\rr_r\log \rr_r - \rr_t\log \rr_t \big\}(y)\d y+ 2K^2(2K^2+1)(t-r).\end{split}\end{equation}

(b) In general, let $0\le \psi\in C_0^\infty(\R^d)$ such that $\int_{\R^d}\psi(x)\d x=1$, and define the smooth mollifier $\scr S_n$:
$$\scr S_n f(x):= n^d \int_{\R^d} f(x-y) \psi(ny)\d y,\ \ n\ge 1, f\in L^1_{loc}(\R^d).$$
Let
$$b_1^{n)}(t,\cdot):=  \scr S_n b_1(t,\cdot), \ \ \ a_1^{n)}(t,\cdot):= \scr S_n a_1(t,\cdot),\ \ n\ge 1.$$
 Then $(a_1^{n)},b_1^{n)})$ satisfies \eqref{GH}  and  $(A_1)$
 for the same constant $K$.
  So, by step (a) and Lemma \ref{LN1}, the density function $\rr_t^{n)}$ for the diffusion process generated by    $L_t^{a_1^{n)},b_1^{n)}}$ satisfies
\beq\label{YD} \int_r^t\d s\int_{\R^d}\ff{|\nn \rr_s^{n)}|^2} {\rr_s^{n)}} (y)\d y\le c \log \Big(1+\ff t r\Big),\ \ 0<r\le t\le T, n\ge 1\end{equation} for some constant $c=c(K,T,d)>0.$
   Equivalently, for any
$$f\in C_0^{0,2} ([r,t]\times \R^d):=\big\{f\in C_b([r,t]\times\R^d):\  \nn f, \nn^2 f \in C_0([r,t]\times\R^d)\big\},$$
we have
\beg{align*} &   \bigg|\int_{[r,t]\times\R^d} \rr_{s}^{(n)}(y)\DD f_s(y)\d s\d y \bigg|^2 =  \bigg|\int_r^t \d s\int_{\R^d}\big\{\<\nn\log \rr_s^{n)},\nn  f_s\>\rr_s^{n)}\big\} (y)\d y \bigg|^2 \\
&    \le  c \log \Big(1+\ff t r\Big) \int_{[r,t]\times \R^d} |\nn f_s|^2(y)\rr_s^{n)}(y)\d s\d y,\ \ n\ge 1.\end{align*}
By   \cite[Theorem 11.1.4]{SV},
$$\lim_{n\to\infty}   \int_{ \R^d} \rr_{s}^{n)}(y)g(y) \d y=
\int_{\R^d} \rr_{s}(y)g (y) \d y,\ \ g\in  C_b(\R^d),\ \ s\in [r,t].$$
So,  the above estimate   implies
$$\bigg|\int_{[r,t]\times\R^d} \rr_{s} (y)\DD f_s(y)\d s\d y \bigg|^2\le  c  \log \Big(1+\ff t r\Big) \int_{[r,t]\times \R^d} |\nn f_s|^2(y)\rr_s(y)\d s\d y$$ for any $ f\in C_0^{0,2} ([r,t]\times \R^d).$
Therefore, \eqref{YD1} holds.

(c)  If \eqref{XY} holds, then by Malliavin's calculus, see for instance \cite{Nu} or \cite[Remark 2.1]{WZ}, for any $v\in \R^d$ with $|v|=1$, there exists a martingale $M_t^{1,x,v}$ such that
$$\E[\nn_v f(X_t^{1,x})]= \E[f(X_t^{1,x}) M_t^{1,x,v}],\ \ f\in C_b^1(\R^d), t\in (0,T]$$ and
$\E[|M_t^{1,x,v}|^2]\le \ff {c} t$ holds for some constant $c=c(T,K,d)>0$ and all $t\in (0,T].$
This implies
$$\bigg|\int_{\R^d}\big\{ \<v,\nn_x \log p_t^{1,x}\> f\big\}(y) p_t^{1,x}(y)\d y\bigg|^2\le \ff c t \int_{\R^d} f(y)^2 p_t^{1,x}(y)\d y,\ \ f\in C_b^1(\R^d),\ |v|=1.$$
Equivalently,
$$\int_{\R^d}\ff{|\nn p_t^{1,x}|^2}{p_t^{1,x}}(y)\d y\le \ff {c d} t,\ \ \ \ t\in (0,T],$$ so that
\eqref{YD3} holds.
\end{proof}

We are now ready to prove Theorem \ref{T1}.

\beg{proof}[Proof of Theorem \ref{T1}]

(1) Let $p>1$ and $\vv\in (0,\ff 1 2]$ be in Proposition \ref{P2.1}.
By Proposition \ref{P3.1} and $(A_1)$, $(H)$ holds for $\nu=\dd_{x_1}$ and $(a^{\<t_0\>}, b^{\<t_0\>})$ replacing $(a_2,b_2).$
By \eqref{BOG} with $\nu=\dd_{x_1}$ and \eqref{YD1}, we find a constant $c_1=c_1(K,T,d,\varphi)>0$ such that
\beq\label{YD1'} \beg{split}&\Ent(P_{t_1}^{1,x_1}|P_{t_1}^{\<t_0\>x_1})\\
&\le c_1 \bigg[ |a_1-a_2\|_{\vv t_1\to t_1,\infty}^2 + \int_{\vv t_1}^{ t_1}\big(\|{\rm div}(a_1-a_2)\|_{t,\infty}^2
+\|b_1-b_2\|_{t}^2\big)\d t \bigg],\\
&\qquad  \ t_1\in (0,T], x_1\in\R^d.\end{split}\end{equation}
Combining this with
  \eqref{YPP} and Proposition \ref{P2.1}, we find a constant $c=c(K,T,d,\varphi)>0$ such that for any
  $t_1\in (0,T]$ and $x_1,x_2\in\R^d,$
\beg{align*} \Ent(P_{t_1}^{1,x_1}|P_{t_1}^{2,x_2})  \le&\, I_{t_1}(x_1,x_2):=
   \ff {c} {t_1} \bigg(|x_1-x_2|^2 + \int_0^{t_1}\big\{ \|b_1-b_2\|_{s,\infty}^2+\|a_1-a_2\|_{s,\infty}^2\big\}\d s\bigg)\\
& + c\bigg(  \|a_1-a_2\|_{\vv t_1\to t_1,\infty}^2+
   \int_{\vv t_1}^{t_1} \|{\rm div}(a_1-a_2)\|_{s,\infty}^2\d s\bigg).\end{align*}
Equivalently, for any $t\in (0,T]$ and $f\in\B_b^+(\R^d),$
\beq\label{AP0} \int_{\R^d} \big\{\log f (y)\big\} P_t^{1,x_1}(\d y)\le \log \int_{\R^d} f(y) P_t^{2,x_2} (\d y) + I_t(x_1,x_2),\ \ x_1,x_2\in\R^d.\end{equation}
Let $\pi\in\C(\nu_1,\nu_2)$ such that
$$\W_2(\nu_1,\nu_2)^2=\int_{\R^d\times\R^d}|x_1-x_2|^2\pi(\d x_1,\d x_2).$$
we obtain
\beg{align*} &\Ent(P_t^{1,\nu_1}|P_t^{2,\nu_2})=\sup_{0<f\in\B_b(\R^d)}\bigg\{ \int_{\R^d} \big\{\log f (y)\big\} P_t^{1,\nu_1}(\d y)- \log \int_{\R^d} f(y) P_t^{2,\nu_2} (\d y)\bigg\}\\
 &\le  \int_{\R^d\times\R^d}I_t(x_1,x_2)\pi(\d x_1,\d x_2)\\
 &= \ff c t\bigg(
  \W_2(\nu_1,\nu_2)^2+\int_0^t \big\{\|b_1-b_2\|_{s,\infty}^2+\|a_1-a_2\|_{s,\infty}^2\big\}\d s\bigg)\\
  &\qquad  +
  c\bigg(  \|a_1-a_2\|_{\vv t\to t,\infty}^2+
   \int_{\vv t}^t\|{\rm div}(a_1-a_2)\|_{s,\infty}^2\d s\bigg).\end{align*}
   Hence, \eqref{AE1} holds.

(2)  Let  \eqref{XY} hold.    By \eqref{BOG} and \eqref{YD3},    we find a constant $c_1=c_1(K,T,d,\varphi)>0$ such that  for any  $t\in (0,T]$ and $x_1\in\R^d,$
\beg{align*}&\Ent(P_{t}^{1,x_1}|P_{t}^{\<t_0\>x_1})
\le c_1  \int_{\vv t}^{t} \ff 1 {s} \|a_1-a_2\|_{s,\infty}^2\d s + c_1\int_{\vv t}^{t}\big[\|{\rm div}(a_1-a_2)\|_{s,\infty}^2
+ \|b_1-b_2\|_{s,\infty}^2\big]\d s,\\
&\le \ff{c_1}{\vv t} \int_{\vv t}^{t} \|a_1-a_2\|_{s,\infty}^2\d s + c_1\int_{\vv t}^{t}\big[\|{\rm div}(a_1-a_2)\|_{s,\infty}^2
+ \|b_1-b_2\|_{s,\infty}^2\big]\d s.\end{align*}
Then as explained above that using this estimate to replace \eqref{YD1'}, we derive \eqref{AE3} for some constant $c=c(K,T,d,\varphi)>0.$\end{proof}

\section{Application to McKean-Vlasov  SDEs }

As an application of  Theorem \ref{T1}, we are able to establish   \eqref{EC}  for \eqref{E1} with distribution dependent multiplicative noise.
For any $\mu\in C([0,T];\scr P_2)$, let
$$a^\mu(t,x):= \ff 1 2 (\si\si^*)(t,x,\mu_t),\ \ \ b^\mu(t,x):= b(t,x,\mu_t),\ \ (t,x)\in [0,T]\times\R^d.$$
Correspondingly to $(A_1)$ and $(A_2)$, we make the following assumption.

\beg{enumerate} \item[$(B)$] There exists a constant $K>0$ such that $a^\mu$ and $b^\mu =b^{\mu,0} +b^{\mu,1} $ satisfy the following conditions.
\item[$(1)$]  For any $\mu\in C([0,T];\scr P_2)$,  $b^\mu$ is locally bounded, and
 for any $ (t,x,\mu)\in [0,T]\times\R^d\times\scr P_2,$
 \beg{align*}  \|\nn b^{\mu,1}\|_{0\to T,\infty}  +\|a^\mu \|_{0\to T,\infty}
 +\|(a^\mu)^{-1}\|_{0\to T,\infty}+\|\nn a^\mu\|_{0\to T,\infty} \le K.\end{align*}
 \item[$(2)$] There exists   $\varphi\in \D$ such that
 $$ \|b^{\mu,0}\|_{T,\varphi} \le K,\ \ \ \mu\in C([0,T];\scr P_2).$$
 \item[$(3)$] For any   $\nu,\mu\in \scr P_2,$
 $$\|b^{\nu} -b^{\mu} \|_{0\to T,\infty}\lor \|a^{\nu} -a^{\mu} \|_{0\to T,\infty}\lor \big\|{\rm div} (a^{\nu} -a^{\mu}) \big\|_{0\to T,\infty} \le K\W_2(\nu,\mu).$$
 \end{enumerate}

\beg{thm}\label{T2}  Assume $(B)$. Then $\eqref{E1}$ is well-posed for distributions in $\scr P_2$, and there exists a constant $c=c(K,T,d,\varphi)>0$ such that
\eqref{EC} holds.
\end{thm}

 \beg{proof}
By $(B)$, for any $\mu\in \scr P_2$, $b^\mu(t,x):=b(t,x,\mu)$ has decomposition $b^{0,\mu}+b^{1,\mu}$ such that $b^{1,\mu}$ is locally bounded and
$$|b^{0,\mu}|\lor \|\nn b^{1,\mu}\|\le K.$$
Let $b^{(1)}:= b^{1,\dd_0},$ where $\dd_0$ is the Dirac measure at $0$,  and let
$b^{(0,\mu)}:= b^\mu-b^{(1)}.$ Then $(B)$ implies
$$|\nn b^{(1)}|\le K,\ \ |b^{(0,\mu)}|\le K+K\mu(|\cdot|^2)^{\ff 1 2}.$$
This together with the the condition on $\si$ included in $(B)$ implies assumptions $(A_0)$ and $(A_1)$ in \cite{HW220} for $k=2$. Therefore, by \cite[Theorem 1.1]{HW220}, \eqref{E1} is well-posed for distributions in $\scr P_2$, and there exists a constant $c>0$ such that
\beq\label{*M} \sup_{t\in [0,T]} \E[|X_t|^2]\le c(1+\E[|X_0|^2])<\infty\end{equation}
holds for any solution with $\L_{X_0}\in \scr P_2.$

For $\nu_i\in \scr P_2, i=1,2,$ and $(t,x)\in [0,T]\times \R^d,$ let
\beq\label{BAI}\beg{split}& a_i(t,x):= a(t,x,P_t^*\nu_i)=\ff 1 2 (\si\si^*)(t,x, P_t^*\nu_i),\\
 & b_i(t,x):= b(t, x, P_t^*\nu_i),\ \ b_i^{(k)}(t,x):= b_i^{k,P_t^*\nu_i}(t,x), \ \ k=0,1. \end{split}\end{equation}
By Theorem \ref{T1}, under $(B)$, there exists a constant $c_1=c_1(K,T,d,\varphi)>0$  such that for any $t\in (0,T],$
\beg{align*} &\Ent(P_t^*\nu_1|P_t^*\nu_2)\le \ff {c_1}t \W_2(\nu_1,\nu_2)^2 \\
&+ c_1\|b_1-b_2\|_{t,\infty}^2+ c_1\log(1+t^{-1}) \|a_1-a_2\|_{t,\infty}^2
+ c_1 t\|{\rm div}(a_1-a_2)\|_{t,\infty}^2\\
&\le  \ff {c_1}t \W_2(\nu_1,\nu_2)^2 + c_1K^2\big\{1+\log(1+t^{-1})+ t\big\}\sup_{s\in [0,t]}\W_2(P_s^*\nu_1,P_s^*\nu_2)^2.\end{align*}
Then there exists a constant $c_2=c_2(K,T,d,\varphi)>0$ such that
$$ \Ent(P_t^*\nu_1|P_t^*\nu_2)\le \ff {c_1}t \W_2(\nu_1,\nu_2)^2 + \ff{c_2} t \sup_{s\in [0,t]}\W_2(P_s^*\nu_1,P_s^*\nu_2)^2,\ \ t\in (0,T].$$
Combining this with the following Proposition \ref{P4.1}, we derive \eqref{EC} for some constant $c>0,$ and hence finish the proof. \end{proof}

\beg{prp}\label{P4.1} Assume $(B)$. Then there exists a constant $c>0$ such that
$$\W_2(P_t^*\nu_1,P_t^*\nu_2)\le c \W_2(\nu_1,\nu_2),\ \ t\in [0,T], \nu_1,\nu_2\in \scr P_2.$$
\end{prp}
\beg{proof} Let $a_i$ and $b_i$ be in \eqref{BAI}, and let $u_t$ be in \eqref{W3} for large enough $\ll>0$ such that \eqref{W4} holds.
Let $X_0^1,X_0^2$ be $\F_0$-measurable such that
\beq\label{R1} \L_{X_0^i}=\nu_i,\ \ i=1,2,\ \ \E[|X_0-X_0^2|^2]=\W_2(\nu_1,\nu_2)^2.\end{equation}
Let $X_t^i$ solve \eqref{SDE} with initial value $X_0^i$. We have $\L_{X_t^i}=P_t^*\nu_i$, so that
\beq\label{R2} \W_2(P_t^*\nu_1,P_t^*\nu_2)^2\le \E[|X_t^1-X_t^2|^2],\ \ t\in [0,T].\end{equation}
Let $\tt X_t^i=X_t^i+u_t(X_t^i),i=1,2.$ Then
\beq\label{W5'} \ff 1 2 |X_t^1-X_t^2|\le |\tt X_t^1-\tt X_t^2|\le 2 |X_t^1-X_t^2|, \ \ t\in [0,T],\end{equation} and similarly to \eqref{SV1}, by \eqref{W3}, \eqref{E1} for $X_t^i$ and It\^o's formula, we have
  \beg{align*}&\d \tt X_t^{1}= \big\{\ll u_t+b_1^{(1)}(t,\cdot)\big\}(X_t^{1})\d t
 + \big\{I_d +\nn u_t(X_t^{1})\big\}\si_1(t,X_t^{1})\d W_t,\\
 &\d \tt X_t^{2}= \big\{\ll u_t+ (L_t^{a_2,b_2}-L_t^{a_1,b_1})u_t + (b_2-b_1^{(0)})(t,\cdot)\big\}
(X_t^{2})\d t\\
&\qquad\qquad\qquad\qquad
 + \big\{I_d +\nn u_t(X_t^{2})\big\}\si_2(t,X_t^{2})\d W_t. \end{align*}
 Combining this with $(B)(1)$, \eqref{W4}, \eqref{R1} and It\^o's formula,  we find $k_1=k_1(K,T,d,\varphi)>0$ such that
$$ \d |\tt X_{t}^{1}-\tt X_{t}^{2}|^2\le k_1\big(|\tt X_{t}^{1}-\tt X_{t}^{2}|^2+\|a_1-a_2\|_{t,\infty}^2+\|b_1-b_2\|_{t,\infty}^2\big)\d t + \d M_t,\ \ t\in [0,T].$$
Noting that $(B)(3)$ and \eqref{BAI} imply
$$\|a_1-a_2\|_{t,\infty}^2+ \|b_1-b_2\|_{t,\infty}^2\le 2 K^2 \xi_t,\ \ \xi_t:=\sup_{s\in [0,t]}\W_2(P_s^*\nu_1,P_s^*\nu_2)^2,$$
and  due to  \eqref{W4}, \eqref{R1} and \eqref{R2}
$$\E[|\tt X_0^1-\tt X_0^2|^2]\le 4 \W_2(\nu_1,\nu_2)^2,\ \ \E[|\tt X_t^1-\tt X_t^2|^2]\ge \ff 1 4
\E[|X_t^1-X_t^2|^2]\ge \ff 1 4 \W_2(P_t^*\nu_1,P_t^*\nu_2)^2,$$
we find a constant $k_2=k_2(K,T,d,\varphi)>0$ such that
$$\xi_t\le k_2 \W_2(\nu_1,\nu_2)^2+k_2 \int_0^t\xi_s\d s,\ \ t\in [0,T].$$
Since \eqref{*M} implies $ \xi_t<\infty,$ by Gronwall's inequality, this implies
$$\sup_{t\in [0,T]} \W_2(P_t^*\nu_1,P_t^*\nu_2)^2=\xi_T\le k_2\e^{k_2T} \W_2(\nu_1,\nu_2)^2.$$
So, the proof is finished.
\end{proof}

\paragraph{Remark 4.1.} After an earlier version of this paper is available online,
the bi-coupling argument developed here has been applied in \cite{HRW, QRW} for singular and degenerate models, where condition $(B)(3)$ is weakened in \cite{HRW} by using $\W_\psi+\W_k$ replacing $\W_2$, see
\cite[Theorem 1.3, Remark 1.2]{HRW} for details. We believe that with additional efforts this new coupling argument will enables one   to derive the entropy-cost inequality for McKean-Vlasov SDEs with singular potentials, where $b_t(x,\mu)$ is given by
$$b_t(x,\mu):= \int_{\R^d} V(x-y)\mu(\d y) $$
for $V$ being a singular potential such as the Coulomb potential
$V (x) =  |x|^{2-d}$ for $d>2$ and $V(x)=\log |x|$ for $d=2.$ This will be addressed in a forthcoming paper.

 \paragraph{Acknowledgement.} The authors  would like to thank Professor Xing Huang and the referees for useful conversations and corrections.

\paragraph{Author Contributions.} The manuscript has been developed in collaboration by the two stated authors.

\paragraph{Data Availibility.} No datasets were generated or analysed during the current study.

\paragraph{Declarations Conflict of interest.} Partial financial support was received from the National Key R\&D Program
of China (NO. 2022YFA1006000, 2020YFA0712900) and NNSFC (11921001).

\end{document}